\documentclass[12pt]{amsart}

\usepackage{graphicx}
\usepackage{cite}
\usepackage{amssymb,amsmath,amsthm}

\theoremstyle{plain}

\newtheorem{theorem}{Theorem}[section]
\newtheorem{proposition}[theorem]{Proposition}

\newtheorem*{maintheorem}{Theorem~\ref{broken windows only}}

\theoremstyle{remark}

\newcommand{\newoperator}[2]{\DeclareMathOperator{#1}{#2}}
\renewcommand{\H}{\qopname\relax o{H}}
\newoperator{\AH}{AH}
\newoperator{\AJ}{AJ}
\newoperator{\BS}{BS}               
\newoperator{\GF}{GF}               
\newoperator{\GH}{GH}
\newoperator{\GJ}{GJ}
\newoperator{\Isom}{Isom}
\newoperator{\KS}{KS}               
\newoperator{\ML}{ML}
\newoperator{\Nbhd}{Nbhd}
\newoperator{\PL}{PL}
\newoperator{\QF}{QF}
\newoperator{\QH}{QH}
\newoperator{\QJ}{QJ}
\newoperator{\QpH}{Q^\prime H}
\newoperator{\Rad}{Rad}
\newoperator{\Tr}{Tr}
\newoperator{\Vol}{Vol}
\newoperator{\area}{area}
\newoperator{\base}{base}
\newoperator{\codim}{codim}
\newoperator{\const}{const}
\newoperator{\degree}{degree}       
\newoperator{\glue}{glue}
\newoperator{\graft}{graft}
\newoperator{\inj}{inj}             
\newoperator{\interior}{int}
\newoperator{\length}{length}
\newoperator{\llength}{\underline length}
\newoperator{\mass}{mass}
\newoperator{\ls}{ls}
\newoperator{\qf}{qf}
\newoperator{\skin}{skin}
\newoperator{\teich}{\textsl{T}}
\newoperator{\unglue}{unglue}
\newoperator{\wb}{wb}               
\newoperator{\window}{window}
\newoperator{\ws}{ws}               
\newoperator{\bom}{\boundary_0(M)}
\newoperator{\dl}{dl}
\newoperator{\dr}{dr}
\newoperator{\Acyl}{Acyl}
\newoperator{\acyl}{Acyl}

\newcommand{\hy}{\mathbb{H}}
\newcommand{\integers}{\mathbb{Z}}
\newcommand{\complexes}{\mathbb{C}}
\newcommand{\proj}{\mathbb{P}}
\newcommand{\reals}{\mathbb{R}}

\newcommand{\refto}[1]{\cite{#1}}
\newcommand{\refin}[2]{#1 of \cite{#2}}
\newcommand{\fund}[1]{\pi_1(#1)}
\newcommand{\set}[1]{\{#1\}}

\newcommand{\Chat}{\hat{\complexes}}

\newcommand{\arrow}{\rightarrow}
\newcommand{\boundary}{\partial}

\newcommand{\compose}{\circ}
\newcommand{\cross}{\times}
\newcommand{\cusps}{{\hbox{\rm cusps}}}

\newcommand{\intersect}{\cap}
\newcommand{\inverse}{^{-1}}

\newcommand{\sinfty}{S_\infty^2}

\newcommand{\tps}{\Chat-\set{0,1,\infty}} 
\newcommand{\union}{\cup}

\begin{document}
\title[Hyperbolic Structures III: incompressible boundary]{Hyperbolic Structures on 3-manifolds, III:\\
Deformations of 3-manifolds with incompressible boundary}
\author{William P. Thurston}
\date{1986 preprint $\to$ 1998 eprint}
\begin{abstract}
This is the third in a series of papers construting hyperbolic structures
on Haken manifolds, analyzing the mixed case of 3-manifolds that
with incompressible boundary that are not acylindrical,
but are also not interval bundles. The main ingredient (beyond
\cite{Thurston:hype1} and \cite{Thurston:hype2} is an upper bound
for the hyperbolic length of the `window frame',
that is, the boundaries of $I$-bundles in the manifold, using a 
growth rate estimate for branched hyperbolic surfaces.

This is slightly revised from the 1986 version of 
a preprint that circulated in the early '80's. A few figures have
been added, and a few clarifications have been made in the text.
\end{abstract}
\maketitle
\setcounter{section}{-1}
\section{Introduction}
\label{Introduction}

In the first two parts of this series, \refto{Thurston:hype1} and \refto{Thurston:hype2}, we analyzed
hyperbolic structures on two opposite classes of $3$-manifolds.  In the first
paper, we studied acylindrical $3$-manifolds. In the second, we studied
$3$-manifolds of the form $S \cross I$: such a manifold is really one big
cylinder.

We will now study the hybrid case, that is, general $3$-manifolds with
incompressible boundary.  We will also obtain some information about hyperbolic
structures on three-manifolds  whose boundary is not incompressible ---
provided we bound the lengths of a set of curves on $\boundary M$ sufficient to
intersect the boundary of any essential disk.

Recall from \refin{\S7}{Thurston:hype1} that a {\it pared manifold} is a pair $(M, P)$,
where $P \subset \boundary M$ is a (possibly empty) $2$-dimensional submanifold
with boundary such that
\begin{description}
\item[(a)] the fundamental group of each of its components injects into the
fundamental group of $M$.
\item[(b)] the fundamental group of each of its components
contains an abelian subgroup with finite index.
\item[(c)]  any cylinder 
$$C: (S^1 \times I, \boundary S^1 \times I) \arrow (M, P)$$
such that $\fund C$ is injective is homotopic {\it rel} boundary to $P$
\item[(d)] $P$ contains every component of $\boundary M$ which has an abelian subgroupof finite index.
\end{description}
The terminology is meant to suggest that certain parts of the skin of $M$
have been pared off to form parabolic cusps in hyperbolic structures for $M$.

Recall also that
$\H(M,P)$ is the set of complete hyperbolic manifolds $M$
together with a homotopy equivalence of $(M, P)$ to $(N, \cusps)$, where
we represent cusps by disjoint horoball neighborhoods. There are three
topologies $\AH(M,P)$, $\GH(M,P)$
and $\QH(M,P)$ on $\H(M,P)$: the algebraic, geometric, and quasi-isometric
or quasiconformal topologies.

A pared manifold is {\it acylindrical} if $\boundary M - P$ is incompressible and if every cylinder
$$C: (S^1 \times I, \boundary S^1 \times I) \arrow (M, \boundary M - P)$$
such that $\fund C$ is injective is homotopic {\it rel} boundary
to $\boundary M$. 
Theorem \refin{7.1}{Thurston:hype1} asserts
that $\AH(M,P)$ is compact when $(M,P)$ is acylindrical.

This is not the case if $(M,P)$ admits essential cylinders. 

In order
to describe what unbounded behaviour there is for $N \in \AH(M,P)$,
we will make use of
the {\it characteristic submanifold} of $(M,P)$, whose theory was
first developed in \refto{Jaco-Shalen:characteristic}
 and \refto{Johannson:characteristic}.   The characteristic submanifold
is defined provided $\boundary_0(M) = \boundary(M) - P$
is incompressible.  The characteristic
submanifold is a disjoint
union of Seifert fiber spaces and interval bundles
over surfaces
in a $3$-manifold.
It is determined up to isotopy.  The ends of the intervals
are on $\boundary M$, and in our case, we will consider only those interval
bundles where the ends of the intervals are contained in
$\boundary_0(M)$.
The characteristic submanifold has the universal property that is
essential, and that any other essential Seifert fiber
space or interval bundle in $M$ with the given boundary conditions
is isotopic into the characteristic
submanifold.  The definition and development
of a sufficient theory to handle embedded interval bundles and Seifert spaces
is fairly elementary, and akin to the theory of the prime decomposition
of a three-manifold.  The theory also handles homotopy classes of maps
which are not necessarily embeddings; this is more difficult.
We will make free use of both types of universal properties of the
characteristic submanifold, for embeddings and for homotopy
classes, even though we could probably get by with only the theory
for embeddings.

In the case of current interest, that of
a pared manifold which admits at least one hyperbolic structure, the
Seifert fiber space part of the characteristic submanifold consists
of some number of
thickened tori in neighborhoods of torus components
of $P$, and some number of solid tori, each
of which intersects $\boundary_0(M)$ in  a non-empty union of
annuli.  These annuli wind one or more times around the
``long way'' of the solid torus before closing.  A transverse disk
in a solid torus intersects the boundary annuli in two or more arcs.
However, if it intersects in only two arcs, we can and will interpret
the solid torus as an interval bundle, rather than as a Seifert
fiber space.  The transverse disks
in all remaining Seifert fiber spaces intersect $\bom$ in at least
three arcs.

It is convenient for our purposes to carve out of the characteristic
manifold something which we will call the {\it window},
which will be an interval bundle $\window (M,P)$.  If $M$ is made
of glass, then viewed from the outside of $M$
the window is the part that one can see through
without homotopic distortion or entanglement.
To form the window, begin with
all of the interval bundles in the characteristic submanifold.
In addition, a portion of the boundary
of each solid torus in the characteristic submanifold consists of
essential annuli in $(M, \boundary_0(M))$.  Thicken each of these,
to obtain additional interval bundles.  Some of the new interval
bundles, however,
may be redundant: they may be isotopic into other interval bundles.
If such absorbtion is possible, then eliminate interval bundles
which can be absorbed, one by one.  The result is the window.
The ends of the intervals of the window form a subsurface of
$\boundary_0(M)$, the {\it window surface} $\ws(M,P)$; it is a $2$-fold
covering of $\wb(M,P)$, the {\it window base}.  We may assume
that $\window(M,P)$ contains all of the annuli in $P$.
Note that even when $M$ is orientable, the window base need not be
orientable, since non-orientability of the fibration
can cancel with non-orientability of the base to make the total
space orientable.

The main theorem of this paper is
\begin{theorem}[Broken windows only]\label{broken windows only}
If $\Gamma \subset \fund M$ is any subgroup which is conjugate
to the fundamental group of a component of $M - \window(M,P)$,
then the set of representations of $\Gamma$ in $\Isom(\hy^3)$
induced from $\AH(M,P)$ are bounded, up to conjugacy.

Given any sequence $N_i \in \AH(M,P)$, there is a subsurface with
incompressible boundary
$x\subset \wb(M,P)$ and a subsequence $N_{i(j)}$ such that the
restriction of the  associated
sequence of representations $\rho_{i(j)}: \pi_1(M) \to \Isom(\H^3)$
to a subgroup $\Gamma \subset \fund M$ converges
if and only if $\Gamma$ is conjugate to the
fundamental group of a component of $M - X$, where $X$ is the
total space of the interval bundle above $x$.
Furthermore, no subsequence of $\rho_{i(j)}$ converges on any larger
subgroup.
\end{theorem}

Compare with \refin{Theorem 6.2}{Thurston:hype2}.

\begin{figure}[htbp]
\begin{minipage}{.47\textwidth}
\includegraphics[width=\textwidth]{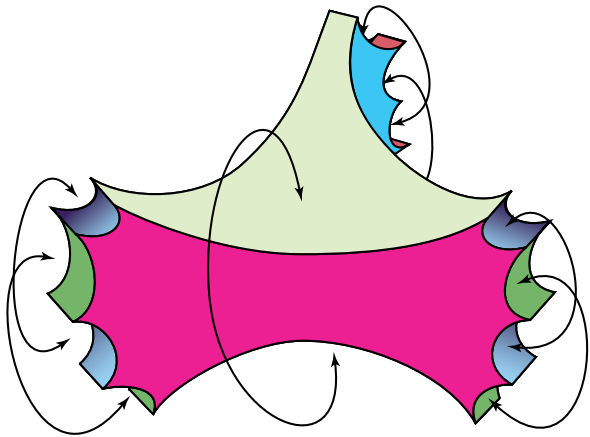}
\end{minipage}
\begin{minipage}{.53\textwidth}
\caption[A ducky 3-manifold]{This is a diagram for
a 3-manifold obtained by gluing
three thickened punctured tori to a solid torus.
Its window consists of the three
thickened punctured tori. Theorem \ref{broken windows only} says that
among all complete hyperbolic structures,
the length of the core circle of the solid torus is bounded. 
Theorem \ref{window frame bounded} and the ensuing discussion
gives an explicit upper bound of 18.13 (which does not pretend to
be near optimal).
}
\label{figure: duck}
\end{minipage}
\end{figure}
The distinction between the characteristic submanifold and the
window is important in this statement.  A good example is a
three manifold obtained by thickening three surfaces with boundary,
and then gluing the annuli coming from the thickened boundary components
to a solid torus (figure \ref{figure: duck}.)
 The characteristic submanifold has a part which
is the solid torus, and three other parts which are the thickened
surfaces.  According to one reasonable definition, the entire manifold
consists of characteristic submanifold (if parallel boundary
faces of the characteristic manifold are pasted together).
Theorem \ref{broken windows only} applied to this example
boils down to the assertion that the length of
the core of the solid torus is bounded in $\AH(M,P)$.

\medskip
The proof divides into two main steps.  The first step is to show that
the total length of geodesics representing $\boundary \wb(M,P)$
in three-manifolds $N \in \AH(M,P)$ is bounded.  This will be
proven in \S\ref{Membrane windows},
by a growth rate argument which actually gives reasonable and
concrete bounds for this total length.

The second main step
is to show that a certain relative deformation
space, that of hyperbolic structures on $M - \window$
relative to $P$ and $\window$, is 
compact.  This is quite parallel to the proof of the
main theorem of part I of this series, \refto{Thurston:hype1}.
It will be a corollary of
a more general theorem we will prove, which asserts that if $(M, \boundary M)$
is a three-manifold with boundary, and if $X \subset \boundary M$
is a sufficiently complicated collection of curves on $\boundary M$, then
when $X$ is held bounded, so is $M$.  See
Theorem \ref{relative boundedness} for the formal statement.
This relative boundedness theorem will make use of the
relative uniform injectivity of pleated surfaces,
proven in \S\ref{Relative injectivity of pleated surfaces}.


\section{Membrane windows}
\label {Membrane windows}

We will construct a certain branched $2$-manifold, which describes
the boundary and how some parts of it are homotopically identified
with other parts.

Branched $2$-manifolds are an analogue of
train tracks, in one dimension higher.
We will make use of them in two forms:
in the form of non-Hausdorff $2$-manifolds,
and also in the form of the Hausdorffications of these $2$-manifolds
which are standard $2$-complexes.

It is easy to construct many examples of non-Hausdorff manifolds:
if you begin with a standard manifold, and then identify
an open set of one by a homeomorphism to an open set in another,
the result is still always a manifold, since each point
has a neighborhood homeomorphic to a subset of
$\reals^n$.  However, the identification space is typically
not Hausdorff, since there are
likely to be pairs of points on the frontiers of the open sets such that
all neighborhoods of one intersect all neighborhoods of the other
in the quotient topology.  Non-Hausdorff manifolds arise frequently
as quotient spaces of dynamical systems.

We will construct a branched manifold $\BS(M,P)$ associated
with a pared manifold having incompressible boundary, which
either is a quotient space of $\boundary_0(M)$, or at least has the
image of
$\boundary_0(M)$ as a dense subset.
This branched manifold admits a hyperbolic structure
such that the branch points lie on a finite number of closed geodesics.
There will be a map $\BS(M,P) \arrow M$ defined up to homotopy, and
the inclusion of $\bom$ into $M$ factors up to homotopy through this
map.

\begin{figure}[htbp]
\begin{minipage}{.3\textwidth}
\includegraphics[width=\textwidth]{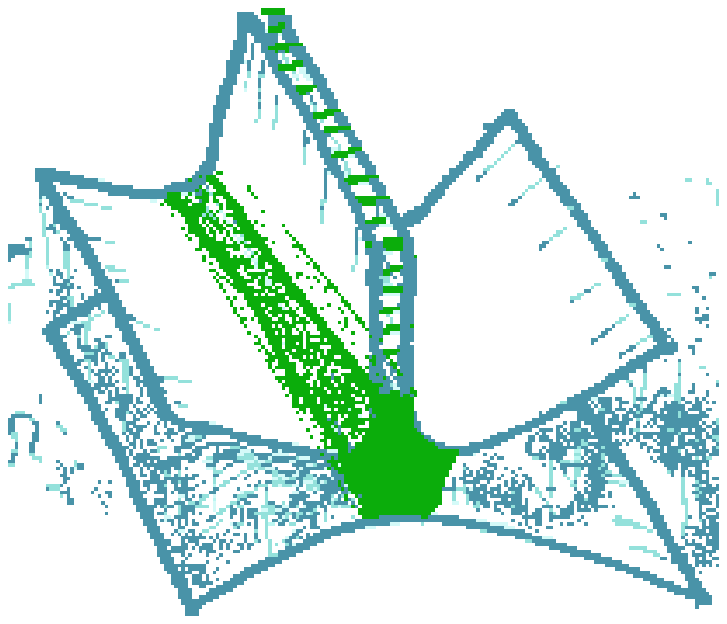}
\includegraphics[width=\textwidth]{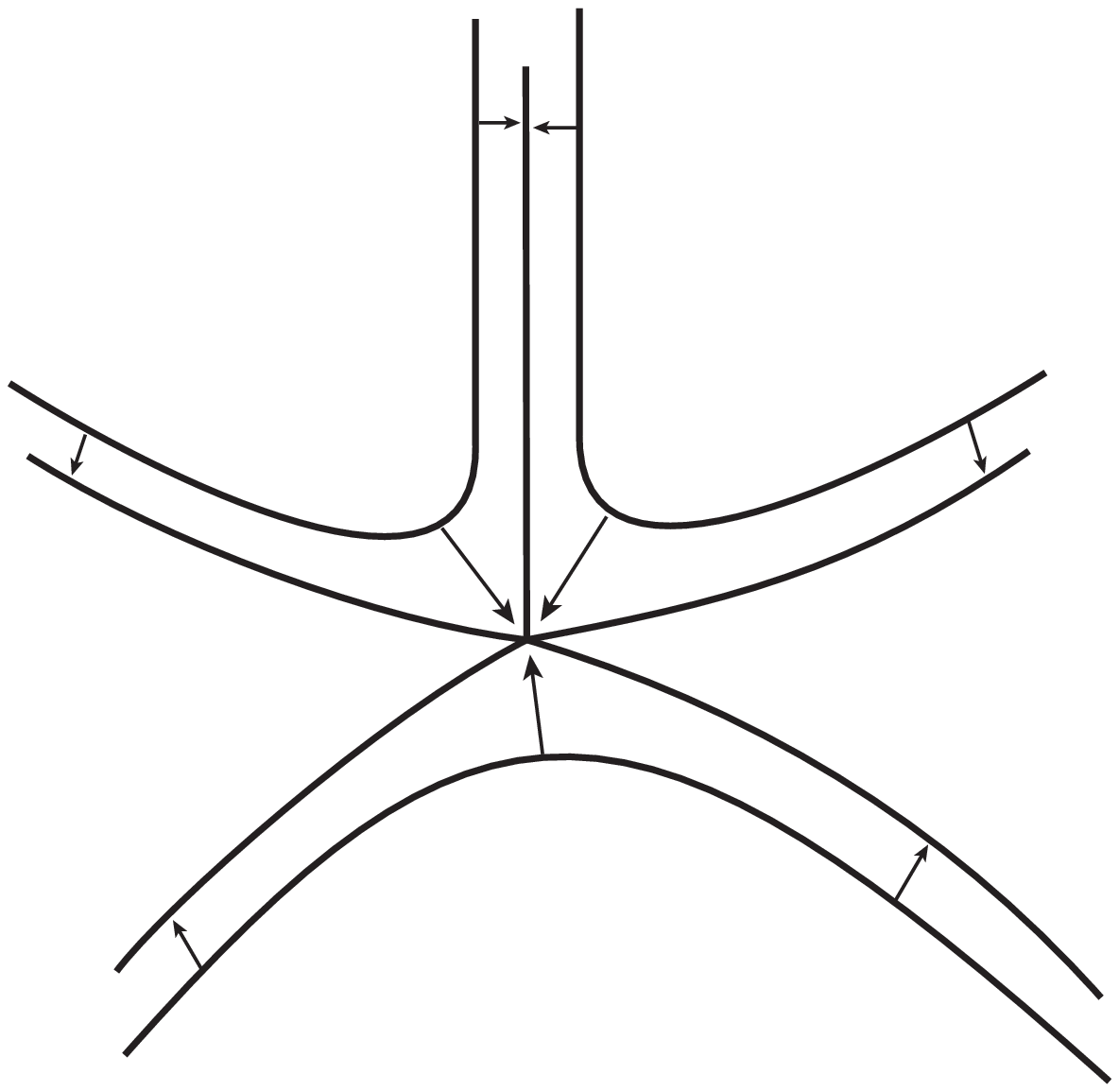}
\end{minipage}
\begin{minipage}{.7\textwidth}
\caption[Constructing KS from windows]{
A certain $2$-complex $KS(M,P)$ is associated with
an atoroidal  pared $3$-manifold $(M,P)$,
roughly described
by starting with $\boundary_0(M)$ and collapsing the $I$-directions of
its windows (characteristic $I$-bundles), and collapsing the
transverse disks of the essential solid torus Seifert fiber spaces.
\emph{Above}, an essential solid torus having three annuli on its
boundary that are part of $\boundary_0(M)$ and three annuli
that cut through the interior of $M$. The uppermost interior annulus
adjoins a `wide' part of $\window(M,P)$, while the other two attach to
portions of $M$ that are relatively acylindrical---the adjoining
windows are just thin slivers of glass.
\emph{Below}, a cross-section of $\KS$, showing the collapsing maps
from $\boundary_0(M)$.
}
\label{figure: constructing KS}
\end{minipage}
\end{figure}

First, we construct the Hausdorff version.
Begin with a characteristic submanifold $W$ for $M$.
Throw out any components of $W$ which are neighborhoods of cusps.
Adjust each solid torus in the
Seifert fiber part of $W$ so that instead of intersecting
$\boundary M$ in a union of annuli, it intersects in a union
of circles.
Similarly, adjust any solid torus which is a component of
the interval bundle part so that any annulus of intersection with $\boundary M$
is collapsed to a circle.
Squeeze together any pairs of annuli in the frontier of
$W$ which are parallel in $M$.  
Now collapse each solid torus to a circle, and collapse 
each of the interval bundles to its base.
The result is a $3$-complex
$K(M,P)$, and the image of $\boundary_0(M)$ in $K(M,P)$
is a $2$-complex $\KS(M,P)$.

$\KS$ is a $2$-manifold except along certain branch curves, some
which
come from the solid tori in $W$ and others
from annuli in the frontiers of the interval bundles.
The components of $\KS$ cut by these curves are surfaces
of negative Euler characteristic. Choose a length for each
branch curve, and extend this choice to give a metric on $\KS$
that is hyperbolic except at the branch curves,
which are geodesics.
\begin{figure}[hbtp]
\begin{minipage}{.33\textwidth}
\includegraphics[width=\textwidth]{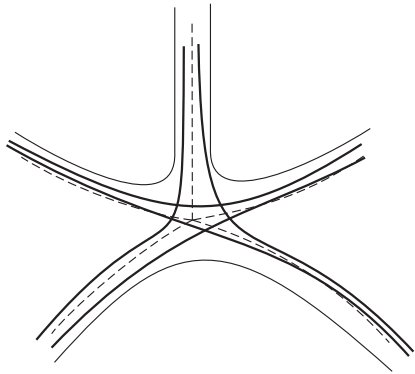}
\end{minipage}
\begin{minipage}{.67\textwidth}
\caption[Construction of branched tunneling boundary surface]{
This is 
the branched surface $\BS(M,P)$ which locally isometrically embeds in
$\KS(M,P)$, from figure \ref{figure: constructing KS}.
Each branch curve $\alpha \subset \KS(M,P)$
is resolved to a number of copies: in this example, there are
3 copies of $\alpha$ coming from the components
$\boundary_0(M)$ that maps to $\alpha$, and 5 copies  for the
pairs of non-adjacent normal directions to $\alpha$.
In a more general example, a cross-section disk could rotate
when it is isotoped once around the solid torus, so that $\alpha$
is replaced by some covering space of $\alpha$,
but this does not change the basic
principle that $\BS(M,P)$ actually branches. }
\label{figure: branched surface}
\end{minipage}
\end{figure}

We will construct a branched surface $\BS$ from $\KS(M,P)$ by replacing
each branch curve $\alpha$ by a probably disconnected
covering space $m(\alpha)$, with a neighborhood on $\BS$ that
locally embed in in $\KS(M,P)$. For each point $t \in \alpha$,
an element in its pre-image is determined by specifying a pair
of sheets of $\KS(M,P)$ at $t$;  these sheets might be
permuted by a rotation when you go once around $\alpha$.
We want to create
as much branching as we can without violating a certain
homotopy injectivity principle to be described below, in
order to get tighter bounds
for the geometry of hyperbolic structures on $(M,P)$,

As a start, we can put the union of all circles on $\bom$ that map to $\alpha$
into $m(\alpha)$.
In cases coming from a window base that is an annulus
or Moebius band, these
boundary circles are not enough to create actual branching. 
But we can add to
$m(\alpha)$, in addition to the boundary circles,
one circle for each pair of sheets of $\KS(M,P)$ along $\alpha$
that are non-adjacent with respect to circular order around $\alpha$.
We define $\BS(M,P)$ by taking the union of the
the complement of the branch circles
in $\KS(M,P)$ with the annular neighborhoods of $m(\alpha)$,
identifying them on open sets. This quotient space $\BS(M,P)$ is a
non-Hausdorff surface, with
a hyperbolic structure inherited from $\KS(M,P)$.

We need to construct
a ``fundamental cycle'' for $\BS$, which can be defined
as a locally finite singular $2$-cycle which has positive degree at
every point. 
A locally finite
$2$-cycle is homologous to one which is transverse to
any particular point, so the degree makes sense even at
the branch curves.
Many non-Hausdorff manifolds without
boundary do not admit such.  For a $1$-dimensional
example, a non-recurrent
train track on a surface determines a non-Hausdorff $1$-manifold
with no fundamental cycle.

In the present case, we can start
with the fundamental cycle of $\bom$, which pushes
forward to a cycle on $\BS$.  It does not have positive degree quite
everywhere:  it has zero degree at each of the circles
which bridge between non-adjacent sheets around a branch curve
$\alpha$ of $\KS$.  A local adjustment can easily be made
so that all these degrees are positive.

Geodesics make perfect sense on a non-Hausdorff hyperbolic surface:
they are paths which have local behavior modeled on a straight
line in the hyperbolic plane.  Unlike in the Hausdorff case,
a tangent vector does not have a unique geodesic extension.
If you construct a geodesic in a certain direction, you may eventually come
to a branch point; at such a place, there may be several
possible extensions of the geodesic.  Our surface
$\BS$ satisfies the property that every geodesic
has an extension to a bi-infinite geodesic; we will call
a non-Hausdorff hyperbolic surface of this type a {\it complete}
surface.

A non-Hausdorff hyperblic surface still has a developing map,
defined by analytic continuation, that maps some covering space
to the hyperbolic plane. 

We can define a space $\GF$ on which the geodesic flow
exists, as the set of pairs $\GF = \set{(X, l)}$, where
$l$ is a bi-infinite unit-speed geodesic on $\BS$, and $X$ is a 
unit tangent vector to $l$.  Give $\GF$ the compact-open topology.
There is a continuous map $p: \GF \arrow T_1(\BS)$, which forgets
$l$.  The preimage of a point $X$ is a totally-disconnected set,
most likely a Cantor set---one can think of $p^{-1}(x)$ as  
the sequence of decisions 
$l$ has to make each time it comes to a choice of branches. 
The space $\GF$ is Hausdorff:  even when
$X_1$ and $X_2$ are tangent vectors along the branch locus
which do not have disjoint neighborhoods,
any two
elements $(X_1, l_1)$ and $(X_2, l_2)$ have neighborhoods
which are disjoint, distinguished from each other in the Cantor set direction.

The map $\GF \arrow T_1(S)$ is not a fibration, because the particular
choices confronting a geodesic which starts out at a tangent
vector $X$ depends on $X$.  For example, some geodesics never
cross branch lines, or cross only a finite number of branch lines,
or cross branch lines in forward time but not in backward time.

A fundamental cycle $Z$ for a branched surface determines an invariant
transverse measure for the geodesic flow,
by the following stipulations:
\begin{description}
\item[(a)] the total transverse measure of the preimage of any point
$p \in \BS$ is
the degree of $Z$ at $p$, and is distributed among the tangent
vectors $X$ in proportion to Lebesgue measure on the circle.
\item[(b)]  Among the elements of $\GF$ over a given tangent vector
$X \in T_1(\BS)$ such that they traverse the
same tangent vectors for time $0 \le t < t_0$ or $0 \ge t > t_0$, 
the probability for the possible choices
of tangent vector $X_{t_0}$ at time $t_0$ is proportional to the degree of
$Z$ at the base point of $X_{t_0}$.
\end{description}

This transverse measure gives rise also to a measure on $\GF$, defined
locally as the product of Lebesgue measure along geodesics with
the transverse measure.
\begin{figure}[htbp]
\centering
\includegraphics[height=1.5in]{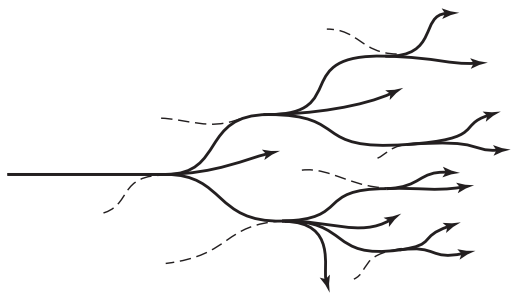}
\caption[Non-deterministic geodesic flow]{Geodesics on non-Hausdorff surfaces
can coincide for a time until they reach a branch point, then separate.
The space $\GF$  is the global description of this behavior,
consisting of all possible  unit tangent vectors to $\BS(M,P)$
together  with a bi-infinite geodesic through the vector.
Proposition \ref{growth proportional to branching} limits the rate
of branching by the growth of volume in $\H^3$. This limit in turn
limits the total length of the branch locus.
}
\end{figure}

\begin{proposition}[Exponential map injective]\label{exponential map injective}
Let $\tilde {\BS}(M,P)$ be the covering space induced from the universal
covering $\tilde M$ of $M$. Let $x \in \tilde{\BS}(M,P)$ be an arbitrary
point.  The exponential ``map'' at $x$ is injective, in the sense that if $l_1$
and $l_2$ are any two geodesic intervals beginning at $x$ which have the same
endpoints, then they coincide.
\end{proposition}

Proposition \ref{acylindrical part doubly incompressible} is a
related statement of a more topological form.

\begin{proof}[Proof of \ref{exponential map injective}]
If there are geodesics $l_1$ and $l_2$
on $\tilde {\BS(M,P)}$ which have identical endpoints, but are not identical,
push
them forward to $\tilde {\KS}(M,P)$.  
Pick a maximal interval of the image of
$l_1$ whose interior does not
intersect the image of $l_2$.  This configuration gives rise
to a rectangle in $ (M - W)$, with two of its sides
on $\boundary M$ and the other two on the characteristic submanifold $W$.
Every such rectangle
can be deformed, rel $W$, into $W$, by the theory of
the characteristic submanifold.
(See \refto {Johannson:characteristic} or \refto {Jaco-Shalen:characteristic}.)
This contradicts the fact
that $l_1$ is a hyperbolic geodesic.
\end{proof}

The locus of endpoints of all geodesics of length $R$ emanating from a point
$x \in \BS$ is a branched $1$-manifold $C_R(x)$.  In an ordinary hyperbolic
surface, this curve is the sphere of radius $R$,
and its total length is $2 \pi \sinh R$.  For
a hyperbolic structure on $\BS$, $C_R(x)$ will tend to grow
faster the more the surface branches.  To make a formal statement,
let $\beta \subset \BS$ be any component of the branch locus.
For any fundamental cycle $Z$, there is a constant value $\degree(\beta)$
to the degree of $Z$ along $\beta$.  Choose an arbitrary
orientation for $\beta$, so we can talk about its left
side and its right side.  The degree of $Z$ just to the
left of $\beta$ is also a constant, $\dl(\beta)$, and
the degree just to the right is another constant $\dr(\beta)$.
\begin{proposition}[Growth proportional to branching]\label{growth proportional to branching}
Let $Z$ as above be a fundamental cycle for $\BS$, and let $g$
be a hyperbolic structure on $\BS$.
Let $Y$ be the set of branch curves on $\BS(M,P)$.
For each $\beta \in Y$, define
$$l(\beta) = -\log \left ( \frac{\degree(\beta)}{\dl(\beta)}\right )
 -\log \left (\frac{\degree(\beta)}{ \dr(\beta)}\right ).$$
Then the average $A(R)$, averaged over $x \in \BS$,
of $\length(C_R(x) \degree_x(Z)$, is at least
$$ 
 2 \pi \sinh(R)
\exp \left ( R 
\sum_{ \beta \in Y}
\frac{.5 \length_g(\beta) \degree(\beta) }{ \mass(Z) }
l(\beta) \right )
$$
\end{proposition}

{\it Remark.} We are really making an estimate for the
entropy of the geodesic flow with respect to the
measure determined by $Z$.  The quantity $l(\beta)$ is minus
the log of the probability
that a geodesic coming along on the right chooses to cross $\beta$,
plus a similar term for the left.  In the
summation, this is multiplied by the
ratio of the flux of the geodesic flow through either side
of $\beta$ to the total measure of $\GF$.

\begin{proof}[Proof of \ref{growth proportional to branching}]
For each $x$ in $\BS$ and for each $R > 0$, we will
define a certain function $F_{(x,R)}$ on the part $\GF_x$
of $\GF$ above $x$, whose integral is the total length
of $C_R(x)$.  $F_{(x,R)}$ will be constant over each set of $(X,l) \in GF_x$
such that the geodesics $l$ agree for time $0 \le t \le R$. These
conditions determines $F$:  the value $F_{(x,R)} ( X,l)$
is $2 \pi \sinh(R) / \degree_x(Z) $ 
times the product, over all choices which $l$ makes
during time $0 \le t \le R$, of the reciprocal of the probability
of making that choice.

The average $A(R)$ of $C_R(x) \degree_x(Z)$, as $x$ ranges over $\BS$, is the
average of $F_{(x,R)} (X,l) \degree_x(Z)$.  The products in the formula for
$F_{(x,R)} (X,l)$ can be thought of as $\exp$ of the sum of minus the
logarithms of the probabilities.  Since $\exp$ is convex upward, $A(R)$ is
greater than or equal to $2 \pi \sinh(R) \exp (B(R))$, where $B(R)$ is the
expected value, among all geodesics of length $R$ in $\BS$, of the sum of minus
the logarithms of the probabilities of the choices it makes. 

There is an exact formula for $B(R)$.  Each branch curve $\beta$ on $\BS$
contributes $-\log(\degree(\beta) / \dl(\beta)$ each time a geodesic flows
through from the left, and the total volume of the flow through $\beta$ in time
$R$ is $ .5 R \length(\beta)\degree(\beta)$. There is a symmetrical formula for
flow from the right. The contribution to $B(R)$ coming from flow through
$\beta$ is therefore
$$\frac{ .5 R \length_g(\beta)}{ \area(\BS)} l(\beta),$$
so 
$$A(R) \ge 2 \pi \sinh(R) \exp \left (
\sum_{\beta} \frac{.5 R \length_g(\beta)\degree(\beta)} {\area(Z)} l(\beta)
\right ).$$
The proposition follows.
\end{proof}

\begin{theorem}[Window frame bounded]\label{window frame bounded}  For any pared manifold
$(M, P)$ such that $\bom$ is incompressible, there is a constant
$C$ such that among all elements $N \in \AH(M,P)$, the length
in $N$ of $\boundary \wb(M,P)$ is less than $C$.
\end{theorem}

\begin{proof}[Proof of \ref{window frame bounded}]
The surface area of a sphere of radius $R$ in $\hy^3$ is
$4 \pi \sinh^2(R)$, so both the volume of the ball
of radius $R$ and the area of the sphere of radius $R$
have exponent of growth $2R$.  (The exponent of growth
of a function $f(R)$
is $\lim \sup \log(f)$; in this case, the $\lim \sup$ is the limit.)

We can represent the surface $\BS$ as a pleated surface in $N$,
by mapping all the branch curves to their corresponding geodesics
(if they have geodesics), and then extending to the remaining
parts of $\BS$, which are ordinary surfaces with boundary.
If some of the branch curves are parabolic elements, so that
they have zero length, then we can similarly construct
a pleated surface based on $\BS$ minus a union of branch geodesics;
it is still a branched surface of a similar type.  This imparts
a complete hyperbolic structure of finite area to $\BS$ or to some subsurface
$Q$.

Fix a fundamental cycle $Z$ on $\BS$, as described.
It defines also a fundamental cycle on any of the subsurfaces $Q$
which might arise.
From proposition \ref{growth proportional to branching},
we see that the exponential growth rate of
$A(R)$ is $1$ plus a positive linear combination of the
lengths of the  curves $\beta$.  (If some branch curves were
removed from $\beta$, then they are not included in the
linear combination, but the inequality still holds since their
length is zero.) 

If $Q$ is compact, the argument is direct.
In that case, there is an upper bound to the area of intersection
of $\tilde Q$ with a ball of radius say $.1$ in $\hy^3$.
The volume of a ball of radius $R$ in $\hy^3$ is
$\int_0^R 4 \pi \sinh^2(t) \,dt$, so it grows with exponential
growth rate $2R$.  Therefore, the image of the exponential
map at any point in $\tilde Q$ can have area at most a constant times
$\exp(2R$), so the total length of the branch set must be bounded.

In general,
each cusp of $Q$
is contained in a subsurface bounded by branch curves, with
no branch curves in the interior.  Construct horoball
neighborhoods of the cusps.  If $x \in \tilde Q$ is any point,
and if $x$ is in a horoball $H$ embedded in $\tilde Q$
which covers
one of the horoball neighborhoods of cusps,
then the area of intersection of the ball of radius $R$ about $x$
with $H$ is less than $2 \pi \sinh(R)$.
Whether or not $x$ is in a horoball,
if $H'$ is one of the horoballs not containing $x$, then the
area of the intersection of the ball of radius $R$ about $x$ with
$H'$ is no more than a constant ($\frac{e}{e - 1}$)
times the area of intersection with
an outer ring of $H'$ of width $1$.

Consider $Q$ minus smaller horoball neighborhoods of the cusps,
shrunk a distance of $1$.
This subset $K$ of $Q$ is compact, so there is a supremum to
the ratio of volume in $N$ with area of intersection with $K$.
Therefore, the length of the intersection of the ball of radius
$R$ about any point in $\tilde Q$ with $K$ is no more than
a bounded multiple of $\exp(2R)$.  By
\ref{growth proportional to branching}, this implies that
the total length of the branch set is bounded.
\end{proof}

A more careful analysis would show that the geodesic flow
for $Q$ is ergodic, and that the
area of the image of the exponential map in
$\tilde Q$ at any two points has a bounded ratio (with bound
depends on the pair of points, though), so that in particular,
the exponential growth rate of area is independent of $x$.

Note that if the hyperbolic structure is turned into a metric
on the Hausdorffication of $Q$, then the growth rate of the
ball of radius $R$ in the metric might be considerably more
than the growth rate of the image of the exponential map.
This is hard to exploit, however, because it is hard to know
what cover of $Q$ is induced from the universal cover of $M$.
The exponential map, in effect, picks out a subset of
the fundamental group of the universal covering
of the fundamental group of $Q$ which
injects no matter which $3$-manifold it arose from.

As an example of how the constants work out in
Theorem \ref{window frame bounded}, consider the three-manifold
of figure \ref{figure: duck}
obtained by gluing three thickened punctured tori to a  solid torus.
In this case, $\BS$ has three branch curves, each of degree $1$,
and the rest of the surface has degree $2$.  For any hyperbolic
structure, the three branch curves have the same length, $a$.
If $\beta$ is any of these curves, then
$l(\beta) = 2 \log(2) \approx  1.38629$.
The exponential growth rate of $A(R)$ is therefore
$1 + 1.5 \length(\beta) l(\beta)/(12 \pi) \approx
1 + .05516 \length(\beta)$.
Consequently, $\length(\beta) < 18.1294$.
If $n$ punctured tori are glued to a solid torus, then
there are $n(n-1)/2$ branch curves in $\BS$, and we can choose
$Z$ to give each of them degree $1/(n-1)$, and the rest of the
surface degree $1$.  Then $l(\beta) = 2 \log(n-1)$, $\mass(Z) = 2 \pi n$,
so the exponential growth rate of $A(R)$ is
$$ 1 +  \frac{ (n-1) \log(n-1)}{8\pi} \length(\beta),$$ so
$\length(\beta) < (8 \pi)/((n-1) \log(n-1))$.
For $n=4$, the length is less than $7.6256$; for $n=7$,
the length is less than $2.3378$.

\medskip
It would be interesting to work out extensions of this growth rate
argument to more general contexts. 

The argument certainly applies directly in one dimension lower,
to train tracks on surfaces, and gives an estimate for the minimum
length of a train track, over all hyperbolic structures on its surface
and all maps of the train track into the surface; this is not
very exciting.  It also applies to faithful representations of
a surface group in an arbitrary homogeneous space.

What could be very useful would be the extension of this analysis to
general incompressible branched surfaces in a $3$-manifold which carries
at least one surface with positive weights.   The problem is that it does
not seem possible to make a hyperbolic pleated surface to
represent a branched surface which
has vertices where branch curves cross.  It can be given a hyperbolic
structure in the complement of such points, but the exponential map
is no longer injective.  Perhaps by estimating the entropy for the
geodesic flow and correcting for the rate of occurrence of multiple
counting because of geodesics which go in opposite ways
around vertices where positive curvature is concentrated,
an estimate could be worked out.

A related problem is to extend the growth rate estimates
in terms of the
length of the branch locus to the case of hyperbolic branched surfaces
with geodesic boundary.  The estimate should involve the length of
the boundary, and degrade as the boundary becomes longer.  In the
case of a Hausdorff hyperbolic surface with geodesic boundary, the
growth rate was studied by Patterson and Sullivan; the exponent of
growth is the same as the Hausdorff dimension of the limit set for the group.


\section{Relative injectivity of pleated surfaces}
\label{Relative injectivity of pleated surfaces}

A $3$-manifold which admits essential cylinders
decomposes into its window, some miscellaneous
solid tori, and an ``acylindrical'' part $\Acyl(M)$.
It is not really acylindrical, however, once the
windows are removed: it is just that 
cylinders in $M$ cannot cross $\Acyl (M)$ in an essential
way.

To express this in a general way, let $N$ be a $3$-manifold,
$f: S \arrow N$ a compact surface mapped into $N$,
and $X \subset S$ a system of non-trivial
and homotopically distinct simple closed
curves on $S$ including all components of $\boundary(S)$.
Then $(S, X, f)$
is {\it incompressible} in $N$ if 
\begin{description}
\item[(a)]
there is no compressing disk
for $(S, f)$ whose boundary is a curve on $S$
which intersects $X$ in one or fewer points.
\end{description}

The triple $(S, X, f)$ is {\it doubly incompressible} in $N$
if in addition
\begin{description}
\item[(b)]
there are no essential cylinders with boundary in
$S - X$,
\item[(c)]
there
is no compressing disk for $(S, f)$ whose boundary is a curve on $S$
which intersects $X$ in
two or fewer points, and
\item[(d)]
Each maximal abelian subgroup of $\fund {S-X}$
is mapped to a maximal abelian subgroup of $\fund N$.\end{description}

There is also a weaker form of (d); $(S, X, f)$ is
{\it weakly doubly incompressible}
if it satisfies (a), (b), (c), and
\begin{description}
\item[(d1)]
Each maximal cyclic subgroup of $\fund {S-X}$
is mapped to a maximal cyclic subgroup of $\fund N$.
\end{description}

This weakening of (d) allows for simple closed curves on $S$
to be homotopic to a $(\integers + \integers)$-cusp of $N$.
The significance of this is that in geometric
limits, cyclic subgroups can turn into $\integers + \integers$.

\begin{proposition}[Acylindrical part doubly incompressible]\label{acylindrical part doubly incompressible}
The triple
$$\left (\Acyl(M),\  
\boundary \wb(M) \intersect \Acyl(M) ,\  
\subset
\right )$$
is doubly incompressible.
\end{proposition}

Compare \ref{exponential map injective}.

\begin{proof}[Proof of \ref{acylindrical part doubly incompressible}]
This is part of the basic theory of the characteristic submanifold.
See \refto{Johannson:characteristic} or \refto{Jaco-Shalen:characteristic}.
\end{proof}

We will extend the results
of \refin{\S5}{Thurston:hype1}
to apply to doubly incompressible triples $(S, X, f)$.

\begin{theorem}[Relative injectivity]\label{relative injectivity}
Let $N$ be a hyperbolic $3$-manifold, and
$(S, X, f)$ weakly doubly incompressible
in $N$.  If $\lambda$ is any maximal lamination
on $S$ containing all curves in $X$ as leaves,
and if $f_\lambda: P_\lambda \arrow N$
is a $\lambda$-pleated surface, then $\lambda$
injects into $\proj(N)$.

This remains true if degeneracies of the pleated
surface are allowed, where closed curves of $\lambda$
map to cusps of $N$.

If the recurrent part of $\lambda$
consists only of closed curves, then at least a degenerate
pleated surface always exists which represents $S \arrow^f X$.
\end{theorem}

\begin{proof}[Proof of \ref{relative injectivity}]
Let $\rho$ be the recurrent part of $\lambda$.
According to \refin{5.5}{Thurston:hype1},
the map $\rho \arrow \proj(N)$, restricted to
non-degenerate leaves,
is a covering map to
its image, and it extends to a map on a small neighborhood
of $\rho$ in $S$ which is a  covering map to its image.

It follows easily from the weak double incompressibility of $(S,X)$
that this can only be the
trivial covering, so that at least $\rho$ embeds in $\proj(N)$:
condition (d1) guarantees that
a closed leaf can map only as a trivial covering to its
image, condition (b) guarantees that no two closed leaves
are mapped to a single closed leaf in the image, and condition
(b) also prevents  components of $\rho$ with more complicated
neighborhoods to map by non-trivial coverings.

Each end of each leaf $l$ of $\lambda - \rho$ is asymptotic with
either one or two leaves of $\rho$ --- one if it is asymptotic
with a closed leaf, two if it is asymptotic with a non-compact
leaf.   We will refer to these two types
of ends as type 1 ends and type 2 ends.

Suppose $l_1$ has the same image in
$\proj(N)$ as $l_2$.  Since $\rho$ embeds,
it follows that  $l_1$ and $l_2$ are asymptotic (on $S$)
at both ends.

A closed loop can be formed, by bridging
between the two leaves at their two ends along
short arcs.  If an end is of type 2, then the short
arc will not to intersect any closed leaves
of $\rho$.  If an end is of type 1, then the arc can be
chosen so that it intersects a closed leaf of $\rho$
in at most one point.  In particular, the closed
loop on $S$ intersects $X$ in at most two points.
It follows from condition (c) that the image of
the loop is null-homotopic in $N$.  It follows
from hyperbolic geometry that $l_1 = l_2$.

There remains still the possibility that a 
leaf $l$ of $\lambda-\rho$ could have image in $\proj(N)$
which is a circle.  Such a situation would force
both ends of $l$ to be of type 2, since the set
of identifications under $\lambda \arrow \proj(N)$
is closed.  This also forces the image of $l$
to coincide with the image of the circle at either
end --- which must be a single circle --- and it
forces the closed leaf to be non-degenerate.
Construct
a loop on $S$ which follows along $l$, crosses a short
bridge to the closed leaf, then unwinds on the closed
leaf, finally crossing a bridge to the other end of $l$.
Since the amount of unwinding is adjustable, it can be chosen
so that the entire configuration maps to a null-homotopic
curve in $N$.  This contradicts 
incompressibility, (a), since a small homotopy of the loop
makes it intersect a closed leaf of $\lambda$ in at most
one point.

To prove that a possibly degenerate
pleated surface always exists provided $\rho$ consists
only of closed leaves, we must show that
each leaf $l$ of $\lambda - \rho$, when pushed forward to
$N$ by a continuous map, can be straightened
out to a geodesic without changing the asymptotic
behavior of its ends. In $\hy^3$, the closed
leaves at either end of $l$ either are covered
by geodesics, or they map to cusps.  If one of the curves is
parabolic and one hyperbolic, the endpoints are automatically
distinct, so $l$ can be straightened. If the
end curves are both hyperbolic, then it follows from
incompressibility of $(S,X)$, as above, that the
endpoints are distinct.

Finally, suppose that $l$ is asymptotic to parabolic
curves at both ends.  Map $l$ to $N$ and lift to $\hy^3$.
If the two endpoints are the same, then $l$ can be
retracted into an arbitrarily small neighborhood of a cusp.
Form a subsurface of $S$ from a neighborhood of $l$ together
with the circle or circles it is asymptotic to.
Choose two non-commuting elements of the fundamental
group of this surface which are represented by loops, such
that one loop (say, parallel to one of the closed leaves)
does not intersect $X$, and the other loop intersects
$X$ in at most one point. (Make it from $l$, with
a short bridge between its ends if they are asymptotic,
or a bridge to the second closed leaf, a traversal
of this leaf, a bridge back to $l$, and a return journey along $l$
back to the basepoint, otherwise.)  The commutator
of the two loops is null-homotopic in $N$;
it intersects $X$ in at most two points,
violating (c).
\end{proof}

\begin{theorem}[Relative uniform injectivity]\label{relative uniform injectivity}
Let $S$ be a compact surface, and  $X$ a
collection of non-trivial, homotopically distinct
simple closed curves on $S$ which includes all boundary components.
Let $B$ and $\epsilon_0$ be positive constants.
Among all pleated surfaces $f: S \arrow N$
($N$ a variable hyperbolic manifold)
pleated along laminations $\lambda$ containing $X$,
where $(S,X,f)$ is doubly incompressible
and the total length of $X$ in $N$ is less than $B$,
the associated maps
$$g: \lambda \arrow \proj(N)$$
are uniformly injective on the $\epsilon_0$-thick part of $S$.
That is,
for every $\epsilon>0$ there is a $\delta>0$ such that for any such $S$,
$N$, $\lambda$ and $f$ and for any two points $x$ and $y \in \lambda$
whose injectivity radii are greater than $\epsilon_0$, if $d(x,y) \ge \epsilon$
then $d(g(x), g(y)) \ge \delta$.

The same statement holds true when degenerate pleated surfaces are
allowed, in which certain closed leaves of $\lambda$ may be parabolic.
\end{theorem}

The uniform injectivity theorem of \refto{Thurston:hype1}
is the special case of this when $X = \boundary S$, and where
the constant $B = 0$ --- that is, where any boundary curves
are parabolic.

\begin{proof}[Proof of \ref{relative uniform injectivity}]
As in \refto{Thurston:hype1}, the main step of the proof is to establish that
the geometric limits of surfaces
which satisfy the hypotheses
are at least weakly doubly incompressible.

Consider a sequence $f_i: S \arrow N_i$ of pleated surfaces,
pleated along $\lambda_i \supset
 X$, such that $(S, X, f_i)$
are doubly incompressible, and for
which the total length of the curves of $X$ in $N_i$
is less than $B$.  Let $g_i$ be the metric induced on $S$ by $f_i$.

For each $i$, let $B_i$ be a collection of points, one in
each component of the $\epsilon_0$-thin part of $S$ with respect
to $g_i$.
and let $E_i$ be a collection of orthonormal frames
at the elements of $B_i$.

There is some subsequence such that sequence of hyperbolic surfaces
$S_i$ defined by the metric $g_i$ on $S$ with respect
to the collection of base frames $E_i$
converges to a geometric limit.
By this we mean that, first, if the universal cover of $S$ is developed into
$\hy^2$, where any of the frames $e_i \in E_i$ is sent
to some fixed base frame in $hy^2$,
then the sequence of images of $\fund S$ converge in the Hausdorff
topology for $\Isom hy^2$: this is the geometric limit from the point
of view of $e_i$.
Second, the Hausdorff limit of the image of
$E_i$ from the point of view of $e_i$
must exist, that is, any other frames which stay within a bounded
distance of $e_i$ should converge.  When these conditions are met,
there is a well-defined, possibly disconnected,
geometric limit with respect to $E_i$.

We can pass to a further subsequence so that the system of
geodesics $X$ and the
laminations $\lambda_i$ also converges geometrically, that is,
in the Hausdorff topology when they are developed onto $\hy^2$
using frames in $E_i$ to get started.

Let $R'$ be the topological surface which is the geometric limit.
On $R'$, denote
the limiting hyperbolic metric $g'$, curve system $Y'$
and lamination $\lambda'$.

Each component of $S-X$ is incompressible in $N_i$.   Since the boundary
components of each such component have bounded length with respect to
$g_i$, there is a non-elementary subgroup of $\fund S-X$ based
at any point $x \in B_i$ generated by loops through $x$ of bounded length.
Therefore, the injectivity radius of $N_i$ at the image of $x$ is
bounded abvoe zero.  The image of each $e_i \in E_i$ in $N_i$ can be extended
uniquely to an orthonormal frame $f_i$ in $N_i$, and there is a further
subsequence so that
the $N_i$ with  collection of base frames $\set f_i $ converges geometrically,
to a hyperbolic manifold $N$
and
so that the maps $f_i: S \arrow N_i$ also converge geometrically
to a $\lambda'$-pleated surface $f': R' \arrow N$.

It may be that there are pairs of cusps of $R'$ which came
from the two sides of a sequence of geodesics on $S$ which
either were degenerate, or grew shorter and shorter in the
sequence.  Form a new topological surface $R$ by joining the two
cusps together along a simple closed curve; we can think of $f$
as a degenerate pleated surface for $R$, if we extend
$\lambda$ by adding a closed leaf for each parabolic curve we adjoin,
and spinning the ends of leaves which tend to these cusps
around the new closed leaves.  Label the new curves as belonging
to $Y$ according to whether they were limits of elements of $X$,
possibly after passing to a subsequence.  Pick a homotopy class $f$ of
maps of the new surface $R$ into $N$, which agrees with the
previous homotopy class $f'$ on  $R'$.  For each degenerate curve, this involves
a free choice of a power of the Dehn twist about the curve, and if
the cusp is a $\integers + \integers$ cusp, the choice of how the
annulus wraps around the torus.

We verify the conditions for double incompressibility of
the limiting pleated surface, with curve system $Y$.

Condition (a) is contained in condition \item[(c)] that there are no compressing
disks which intersect $Y$ in two or fewer points.  Suppose there were
an essential disk for $f:R\arrow N$ with boundary a curve on $R$
meeting $Y$ in two or fewer points.
The map $f$ is approximated by a map of $R$ to the $N_i$, for sufficiently
high $i$.  The approximation at least restricted
to $R'$ factors (up to homotopy) through an embedding $j_i$
of $R$ as a subsurface of $S$ with incompressible boundary.  This factorization
extends over the degenerate curves of $R$ as well, since the
annulus of $R$ and the annulus of $S$ both map into the thin
set of $N_i$; the fundamental group of this component of the thin set
is cyclic, and generated by the core curve of the annulus, so the
two annuli must be homotopic {\it rel} boundary.
The compressing disk would be approximated by a compressing disk for $S$
whose boundary intersects $X$ in only two points, contradicting the
hypothesis that $(S, X, f_i)$ is doubly incompressible.

Condition (b) says that there should be no essential annulus for $f:R \arrow N$
whose boundary is disjoint from $X$.  Any essential annulus would
be approximable by an annulus for $f_i: S \arrow N$, for $i$ sufficiently
high.  We need to check that the approximating annulus is essential, that is,
that its two boundary components are not homotopic on $S$.
If the length of either of the boundary curves on $R$ is greater
than zero, then that boundary curve can be represented by a geodesic
on $R$.  Two geodesics are isotopic on $R$ if and only if they coincide;
the same criterion carries over to the approximating surfaces, so an
essential annulus which has a boundary curve of positive length
would carry over to an essential annulus for the approximations,
contradicting the hypotheses.  If
there were an essential annulus whose boundary components have zero
length on $R$, then when it is carried over to
the approximations, it would still essential
for otherwise we would have added a degenerate curve
joining the two cusps in question when we formed $R$.
Again, this contradicts the hypothesis that $(S, X, f_i)$ is
doubly incompressible.

Condition (d1) also persists in a geometric limit.  Let $H$ be a maximal
cyclic subgroup of the fundamental group of $R-Y$.  The embedding $j_i$
of $R$ as
an incompressible subsurface of $S$ for high $i$ carries $H$ to a maximal
cyclic subgroup of $\fund{S-X}$.
Therefore, its image is a maximal cyclic subgroup
of $\fund N_i$.
If $\alpha$ is any element of $\fund N$
such that some power is a cyclic generator of the image of $H$ in $\fund N$,
then $\alpha$ is
approximated in $N_i$ by an element with a power approximately, and
therefore exactly, equal to
the cyclic generator of the image of $H$ in $N_i$; the power can only be $1$.

We have established that any sequence of pleated surfaces satisfying the
hypotheses of the theorem, has a subsequence converging
to a weakly doubly incompressible pleated surface in a hyperbolic
$3$-manifold.  Theorem \ref{relative uniform injectivity}
follows now from
\ref{relative injectivity}.
\end{proof}


\section{Relative boundness for AH(M)}
\label{Relative boundness for AH(M)}

The relative uniform injectivity theorem of the last section has
a direct application to the boundness of
deformations of a hyperbolic $3$-manifold when the total length of
a sufficiently complicated system of curves on its boundary
is held bounded.

\begin{theorem}[Relative boundedness]\label{relative boundedness}
Let $M$ be a $3$-manifold, and $X$ a collection of non-trivial,
homotopically distinct curves on $\boundary M$ such that
$(\boundary M, X, \subset )$ is doubly incompressible.
Then for any constant $A > 0$, the  subset of $AH(M)$
such that the total length of $X$ not exceeding $A$ is
compact.
\end{theorem}

Note that by setting $X = \emptyset$ we obtain
\refin{1.2}{Thurston:hype1}.
By letting
$X$ be the set of core curves of the parabolic annuli of a pared manifold
and setting $A = 0$
we obtain \refin{7.1}{Thurston:hype1}.

\begin{proof}[Proof of \ref{relative boundedness}]
The proof is much the same as the proof of the main theorem of
\refto{Thurston:hype1}, but we will go through the details for the sake of completeness.  

We may as well assume that $\boundary M$ is non-empty.
Let $Z \supset X$ be a collection of curves which
has at least one element on each boundary component of $M$.
Choose a triangulation $\tau$ of $M$, with one vertex on each element
of $Z$ and such that
each element of $Z$ is formed by one edge of $\tau \intersect \boundary M$.
Choose an orientation for each element of $Z$.

Denote the closed unit ball $\hy^3 \union \sinfty = \overline{\hy^3}$.
For any element $N$ of $AH(M)$, an ideal simplicial map $f_N: M \arrow N$
can be defined in the standard way:
let $\tilde f_N: \tilde M \arrow \overline {\hy^3} \supset \tilde N$
by  map a vertex $v$ of $\tilde \tau$
to the positive fixed point at infinity
for the covering transformation which
generates forward motion along the component of $v$ on $\tilde Z$, and
extend $f_N$ to a piecewise-straight map.

There is a canonical factorization of $f_N | \boundary M$
as an ideal simplicial map $i_N$
to a possibly degenerate hyperbolic structure on $\boundary M$, followed
by a possibly degenerate pleated
surface $p_N$.  The hyperbolic structure $g(N)$ on $\boundary M$
is a metric on
$\boundary M$ minus any torus components, and
possibly with certain curves deleted whose two sides become cusps
in the hyperbolic structure.

Associated with $f_N$ is a certain {\it infinity}
subcomplex $\iota_N$ of $\tau$,
consisting of those simplices such that any copy of them in $\tilde M$
is mapped to a single point by $\tilde f_N$. 
The {\it degenerate simplices} for an ideal simplicial
map are those simplices of which
at least one edge is contained in $\iota$. Thus, a degenerate triangle
collapses either to a line, or to a point at infinity.  A degenerate
$3$-simplex collapses to an ideal triangle, a line, or a point.
A $3$-simplex is not called degenerate if it merely flattens to
a quadrilateral, or maps with reversed orientation.

Let $N(i)$ be any infinite sequence of elements of $AH(M)$ such that
the total length of $X$ is bounded by $A$.  We will extract
a convergent subsequence.  For notational simplicity, each
time we pass to a subsequence, we implicitly relabel it
by the full sequence of positive integers.

First, observe that $\iota_{N(i)} = \iota$
does not depend on $i$, for it can be determined homotopically:
a simplex is in $\iota$ if and only if it is homotopic, {\it rel}
vertices, either to a torus boundary component of $M$, or to
a curve in the collection $Z$.

Recall from \refto{GT3M}, \refto{Thurston:3dgt}
or \refto{Thurston:hype1} that
to each edge of each non-degenerate simplex is associated
an {\it edge invariant} in $\tps$.  The edge invariants
of opposite edges of a three-simplex are equal, and the three
values tend to $0$, $1$ and $\infty$ if the shape of
the simplex degenerates.  We may pass to a subsequence
of $N(i)$ such that the edge invariants of each edge of each
nondegenerate $3$-simplex of $\tau$ converge in $\Chat$.

Following \refin{\S3}{Thurston:hype1},
we define a ``bad complex'' $B$ for the sequence $N(i)$,
which is formed as the union of all $1$-simplices (which are bad
only because they do not yield information), all degenerate simplices,
and within each degenerating $3$-simplex, a quadrilateral
attached along the four edges whose edge invariants are tending
to $0$ and $\infty$.  The good submanifold $G$ for the sequence
$N(i)$ is obtained by deleting a regular neighborhood of $B$,
and the ``cutting surface'' $\kappa$ for $\set N(i) $ is defined
as the internal boundary of $G$, $\kappa = \boundary G - \boundary M$.

An ideal triangle in $N(i)$ serves like a base frame: a certain ``view''
of $\tilde N(i)$ is obtained by sending its three vertices to $0$, $1$,
and $\infty$.  If $\alpha$ is any path in $G$ beginning and ending
on a triangle of $\tau$, it defines a sequence of views of $\tilde N(i)$,
which are related by a sequence of isometries of $\hy^3$.  The submanifold
$G$ was defined in such a way that this sequence of isometries converges.
Note in particular that when a path in $G$
enters a degenerating $3$-simplex
through one triangle and exits through the other accessible triangle,
the isometry that
relates the two views converges to the identity.
Therefore, the sequence of representations $\rho(i)$ for $\fund M$
based at any point $p \in G$ converges when restricted to $\fund {(G,p)}$.

As in the earlier paper, we will augment $G$ together with its map to $M$, in
several steps $g_i: G_i \arrow M$, in such a way that the sequence of composed
representations $G_i \arrow M \arrow \Isom(\hy^3)$ continues to converge, until
finally the map to $M$ is a proper map (one which takes boundary to boundary)
of degree $1$.  We start with $G_0$ be $G$, and $g_0$ the inclusion.

The key point is that the surface $\kappa$ is in a certain sense becoming
smaller and smaller in $N(i)$, as the sequence progresses, so that we can think
of it as if it were only $1$-dimensional.  If we proceed far out in the
sequence, we can homotope the map of $\kappa$ to $N(i)$ so that it lies close
to the $1$-skeleton of the image ideal $3$-simplices, so that the area of the
image of $\kappa$ tends to zero, and it geometrically it looks like a
$1$-complex. In fact, we will analyze in terms of a map to a certain graph
$\gamma$ defined by this geometry.

Toward this end, define a subsurface $\kappa_0 \subset \kappa$ to consist of
the intersection of a regular neighborhood of $\iota$ with $\kappa$, together
with a diagonal band on each side of the twisted quadrilateral of $B$.  The
diagonal bands are to be arranged to ``indicate'' vertices of ideal simplices
which are converging together, from the point of view of its component of $G
\cap$ the simplex; the two diagonal bands on opposite sides of a twisted
quadrilateral of $B$ run in different directions. Each component of $\kappa_0$
will be collapsed to a point to form the vertices of $\gamma$.

The intersection of $\kappa - \kappa_0$ with any $3$-simplex is a rectangle,
with two opposite sides on $2$-faces of the simplex. Form a foliation $F$
transverse of $\kappa - \kappa_0$ so that $\boundary \kappa_0 - \boundary M$
consists of leaves of $F$, and the leaves of $F$ are transverse to the
non-degenerate triangles of $\tau$, by foliating each rectangle with a product
foliation. The leaf space of any rectangle is homeomorphic to the unit interval
$[0,1]$.  The parametrizations can be chosen consistently, so that the leaf
space of a rectangle is identified by an isometry to the leaf space of an
adjoining rectangle.  When this is done, the leaves of $F$ are all compact,
since a leaf enters a rectangle at most once.

The graph $\gamma$ is defined by collapsing each component of $\kappa_0$
to a point, and each leaf of $F$ to a point.  Let $p:  \kappa \arrow \gamma$
be the collapsing map.

Proposition \refin{3.7}{Thurston:hype1} shows (by an easy
argument) that the fundamental group of $p \inverse $ of any edge or vertex of
$\gamma$ has image in $\fund M$ which is abelian.  Define $\gamma_0$ to consist
of those cells $\beta$ of $\gamma$ such whose associated group (the image in
$\fund M$ of $\fund{p\inverse(\beta)} $) is non-trivial. Proposition
\refin{3.8}{Thurston:hype1} concludes that the image of the fundamental group of $p
\inverse$ of any component of $\gamma_0$ in $\fund M$ is abelian: this is
deduced from the nature of commutative subgroups of $\Isom(\hy^3)$.

Let $\gamma_1 = \gamma - \gamma_0$.  An edge $e$ of $\gamma_1$ has one of two
types: either $p \inverse(e)$ is a cylinder, or it is a rectangle. Form the
first enlargement $G_1$ of $G_0 = G$ by attaching a $2$-handle to each cylinder
of the form $p \inverse(e)$, $e$ an edge of $G_1$. The extension $g_1$ of $g_0$
over $G_1 - G_0$ is easy to make over the $2$-handles, since their attaching
curves are null-homotopic in $M$.

The next enlargement $G_2$ of $G_1$ is formed by attaching a
{\it semi-$2$-handle} to each rectangle of the form $p\inverse(e)$, for all
edges $e$ in $\gamma_1$. By definition, a semi-$2$-handles is a copy of $D_+^2
\cross [0,1]$, where $D_+^2$ is the intersection of the unit disk in $\reals^2$
with the closed upper-half-plane. It is attached along the round part of its
boundary to the rectangle in question.  We need to make the extension over the
semi-$2$-handles in such a way that $[-1,1] \cross [0,1]$ maps to $\boundary
M$.

The ends of the leaves of the foliation $F$ in $p\inverse(e)$ are on $\boundary
M$, near two edges of of ideal triangles.  These leaves define a sequence of
changes of viewpoint in $\tilde N(i)$, whose composition converges to a
transformation sending the triangle at one end to the triangle at the other
end, while preserving the edge which is near the leaves. If we use the
factorization of $f_{N(i)} | \boundary M$ as $p_{N(i)} \compose i_{N(i)}$, we
see that the corresponding triangles of $\boundary M - \lambda$ are mapped by
the $\lambda$-pleated map $p_{N(i)} $ close to each other in $N$. Since the
change of view between the two triangles converges as $i \arrow \infty$ and
since the thick part of an ideal triangle of $\boundary M - \lambda$ is always
in the thick part of $\boundary M$, we can choose $\epsilon_0$ small enough
that two points $x_1$ and $x_2$ on these two leaves and in the
$\epsilon_0$-thick part of $\boundary M$ (with respect to $g(N(i))$) are mapped
in $\proj N(i)$ so their distance goes to zero with $i$.

We can now apply Theorem \ref{relative uniform injectivity}, to conclude
that the distance between $x_1$ and $x_2$ goes to zero on $\boundary M$ as
measured by the metric $g(N(i))$.  In particular, for each sufficiently high
$i$, there is a well-defined homotopy class of arcs joining $x_1$ to $x_2$ on
$\boundary M$, and this arc is homotopic to the arc joining $x_1$ to $x_2$ on
$\kappa - \kappa_0$.

Does the homotopy class depend on $i$?  Since the total length of $X$ is
bounded by $A$, no component of $X$ can pass between the leaf of $x_1$ and the
leaf of $x_2$ while they are close together --- these leaves are infinite, and
they remain close for a long time, if $i$ is high.  The arc $\alpha(i)$ joining
$x_1$ with $x_2$ for $N(i)$ does not cross any component of $X$. The loops
$\alpha(j) \union \alpha(i)$ are null-homotopic in $M$, so they must also be
null-homotopic on $\boundary M$, by the incompressibility of $(\boundary M, X,
\subset)$.  Therefore, the homotopy class $\alpha(i)$ is independent of $i$.

Map the semi-$2$-handle to $M$ so that the portion $[-1, 1] \cross [0,1]$ of
its boundary covers a regular neighborhood of $\alpha$, to obtain the extension
$g_2$ of $g_1$ over $G_2$.

Define $\kappa_2$ to be the portion of $\boundary G_2$ which we have not mapped
to $\boundary M$.  There is one component of $\kappa_2$ for each component of
$\gamma_0$  union the vertices of $\gamma$ --- the $2$-handles and
semi-$2$-handles have served to sever $gamma$ along the edges associated to a
trivial group. The image of the fundamental group of any component of
$\kappa_2$ in $\fund M$ is abelian.  Furthermore, $\boundary \kappa_2$ does not
intersect $X$.


We form $G_3$ by attaching three-manifolds to the componenets of $\kappa_2$. A
component $C$ of $\kappa_2$ has one of three types. First, it may be a closed
surface whose fundamental group has image in $\fund M$ isomorphic to $\integers
+ \integers$.  In this case, it is homologous to a cusp, and we can attach a
three-manifold $M_C$ whose boundary is $C \union T^2$, with a map to $M$ having
the same image fundamental group, and realizing this homology. Second, it may
be a closed surface whose fundamental group has image isomorphic to
$\integers$;  then we attach a $3$-manifold $M_C$ whose boundary is $C$, with a
map to $M$ having the same image fundamental group as $C$.  Finally, it may be
a surface with boundary whose image fundamental group is isomorphic to
$\integers$.  From the acylindricity of $(\boundary M, X, \subset)$ it follows
that $C$ is homologous to in $H_2(M, \boundary M)$ to $\boundary M$, by a
$3$-chain represented by a manifold whose image fundamental group is also
$\integers$.  We attach such a manifold $M_C$.

Now we have $G_3$ with a proper map $g_3: G_3 \arrow M$. We will analyze the
geometry of $\boundary g_3: \boundary G_3 \arrow \boundary M$ in conjunction
with the possibly degenerate hyperbolic metric $g_{N(i)}$ on $\boundary M$ and
the pleated map $p_{N(i)}$ to $N(i)$ to show that $\boundary g_3$ and hence
$g_3$ have degree one.

The original good submanifold $G = G_0$ intersects every non-degenerate
triangle of $\tau$ in the complement of a regular neighborhood of its
boundary.  We may isotope so that this regular neighborhood is as small as we
like.  The intersection of $G$ with $\boundary M$ then takes up most of the
area of the the possibly degenerate hyperbolic structure for $\boundary M$.

The boundary of $G_1$ was exactly the same as the boundary of $G_0$.   The
boundary behavior of $G_2$ was modified, on account of the semi-$2$-handles. 
These had the effect of joining sides of ideal triangles in pairs, which by the
relative uniform injectivity theorem were close and nearly parallel.  Only
small bridges were added between sides of triangles, so the map of $\boundary
G_2$ to $\boundary M$ has degree one at most image points, as measured by the
metric $g_{N(i)}$.

The boundary curves of a component $C$ of $\kappa_2$ are formed by a sequence
of segments which cycle through the following types:
\begin{description}
\item[(i)] a long stretch inside an ideal triangle and nearly parallel to one
of its sides,
\item[(ii)] a jump and U-turn to another triangle,
\item[(iii)]  a long stretch in the reverse direction nearly parallel to one of
its sides, and
\item[(iv)] a U-turn near a vertex of the triangle.
\end{description}
These curves can be chosen so that they are geodesics except at the U-turns,
where they have geodesic curvature alternately $\pm \pi$. Therefore, their net
geodesic curvature is small.

If $C$ is of the type which actually intersects $\boundary M$, it lifts to a
cover of $M$ with infinite cyclic fundamental group. Each component of
$\boundary C$ is homologous to a standard representative on $\boundary M$:
contained in a geodesic, if there is a geodesic representing $\integers$,
otherwise contained in one of the curves where the metric degenerates. If we
knew that the boundary components of $C$ were simple curves on $\boundary M$,
or if we knew that they lifted to simple curves in the $\integers$ cover, we
could conclude that each of the curves of $\boundary C$ is homologous to a
standard representative of its class  with a chain of small net area, by an
application of the Gauss-Bonnet theorem.

However, the boundary components of $C$ are not necessarily simple curves, and
the form of the Gauss Bonnet theorem which gives a formula for the area of
chains on a surface with boundary a non-embedded union of curves must involve
information about the amount of topological turning of the boundary curves.

One solution to this difficulty is given in \refin{\S6}{Thurston:hype1}, by constructing a
geometric limit of the gluing patterns for the triangles of $\boundary M -
\lambda$, and observing that the curves of $\boundary C$ are simple curves on
the limit surface. Here is an alternate method:

Each curve $\alpha$ of $\boundary C$ is embedded except in
the short bridges from triangle to triangle.  These bridges cannot
cross the thick part of any other ideal triangle, so they may intersect
other portions of $\alpha$ only when it is near a vertex of some
ideal triangle.  Truncate each ideal triangle in a neighborhood of
each vertex, before constructing $\alpha$, making sure to cut
off a large enough neighborhood so that it does not pass under any
of the bridges.  Now, when $\alpha$ is constructed, it is embedded.

It follows that the part of $\boundary M_C$ which is mapped to
$\boundary M$ has net area near $0$.  The degree of
$g_3 | \boundary G_3$ is an integer $\approx 1$, therefore $= 1$.
Some component of $G_3$ must then map with positive degree.
The sequence of representations restricted to the image of
the fundamental group of
any component of $G_3$ converges, by construction, so it follows
from an application of
\refin{4.3}{Thurston:hype1}
that the sequence of representations of $\fund M$ converges, and hence
that the sequence $N(i)$ converges to a limit in $AH(M)$.
\end{proof}


\section{Proof that only windows break}
\label{Proof that only windows break}

We can now obtain the main theorem by logically assembling
the previous results.  Here is the statement once more:

\begin{maintheorem}
If $\Gamma \subset \fund M$ is any subgroup which is conjugate
to the fundamental group of a component of $M - \window(M,P)$,
then the set of representations of $\Gamma$ in $\Isom(\hy^3)$
induced from $\AH(M,P)$ are bounded, up to conjugacy.

Given any sequence $N_i \in \AH(M,P)$, there is a subsurface with
incompressible boundary
$x\subset \wb(M,P)$ and a subsequence $N_{i(j)}$ such that the
restriction of the  associated
sequence of representations $\rho_{i(j)}: \pi_1(M) \to \Isom(\H^3)$
to a subgroup $\Gamma \subset \fund M$ converges
if and only if $\Gamma$ is conjugate to the
fundamental group of a component of $M - X$, where $X$ is the
total space of the interval bundle above $x$.
Furthermore, no subsequence of $\rho_{i(j)}$ converges on any larger
subgroup.
\end{maintheorem}

\begin{proof}[Proof of \ref{broken windows only}]
Let $M$ be a $3$-manifold with incompressible boundary.
Consider any sequence  $\set {N(i)} $ of elements of $ \AH(M)$
Theorem  \ref{window frame bounded} asserts that the total
length of $\boundary \wb(M,P)$ is uniformly bounded.

Proposition \ref{acylindrical part doubly incompressible}
tells us that  the triple
$$\left (\Acyl(M),
\ \boundary \wb(M) \intersect \Acyl(M),
\ \subset \right )$$
is doubly incompressible.  Theorem \ref{relative boundedness}
implies that the image of the restriction
$$ \AH(M) \arrow \AH(\Acyl(M))$$
has a compact closure.

The analogous restriction map for any solid torus component of $M - \window(M)$
(which is not a part of $\Acyl(M)$)
also has a compact closure, since this is just a question of
the boundedness of the length of a component of $\boundary \wb(M)$.

This proves the first paragraph of the theorem.
For the second paragraph, we apply Theorem
\refin{6.2}{Thurston:hype2} to find a subsequence
such that each for each component $S$ of $\boundary M$, there is a subsurface
$c(S) \subset S$ with incompressible boundary, such that
\begin{description}
\item[(a)] For each component $c(S)_i$ of $c(S)$, the sequence of
representations of its fundamental group converges up to conjugacy.
\item[(b)] If $\Gamma$ is any nontrivial subgroup of $\fund {c(S)}$ such that
for some subsequence of the subsequence, its sequence of representations
converges up to conjugacy, then $\Gamma$ is conjugate to a subgroup of $\fund
{c(S)_i}$ for some $i$.
\end{description}

It follows from (b) that the
subsurface $c(S)$, up to isotopy, includes $S_1 = S - \window(M) \intersect S$.
The window is ``transparent'', so that any closed
curve in $c(S) - S_1$ can be pushed through to the opposite side of the
window, giving a homotopy class of curves on a boundary
component $S'$ of $M$ which must be contained in
$c(S')$.  Also, the sequence of representations of $\fund {M \union c(S)}$
must converge; this intersects the fundamental group of $\boundary M$
based on the opposite side of a window in a large subgroup.  These
constraints, along with
lemma 4.3 of \cite{Thurston:hype1}, yield the
description given in the second paragraph of the statement.
\end{proof}
\bibliographystyle{halpha}
\bibliography{hype}
\end{document}